\documentclass[twoside]{article}

\usepackage{amsfonts}
\usepackage{amsthm}
\usepackage{amsmath}
\usepackage{amssymb}

\title{}
\author{}\date{}

\newtheorem{theorem}{Theorem}[section]
\newtheorem{lemma}[theorem]{Lemma}
\newtheorem{proposition}[theorem]{Proposition}
\newtheorem{remark}[theorem]{Remark}
\newtheorem{definition}[theorem]{Definition}

\numberwithin{equation}{section}

\newcommand{\Om}{\Omega}
\newcommand{\Th}{\Xi}
\newcommand{\om}{\omega}

\newcommand{\te}{\theta}
\newcommand{\vfi}{\varphi}
\newcommand{\be}{\begin{equation}}
\newcommand{\ee}{\end{equation}}

\newcommand{\mcs}{{\mathcal S}}
\newcommand{\mca}{{\mathcal A}}
\newcommand{\mcb}{{\mathcal B}}
\newcommand{\mcl}{{\mathcal L}}
\newcommand{\mct}{{\mathcal T}}
\newcommand{\ha}{\mathcal{H}_{\mathbb{C}}}
\newcommand{\h}{h_{\mathbb{C}}}

\newcommand{\mcn}{\mathcal{N}}
\newcommand{\msn}{\mathcal{S}_\mathcal{N}}

\newcommand{\bdc}{\mathbb{C}}
\newcommand{\bdr}{\mathbb{R}}

\newcommand{\al}{\alpha}

\newcommand{\indn}{_{n\geq 1}}

\newcommand{\bese}{\begin{subequations}}
\newcommand{\ese}{\end{subequations}}
\newcommand{\non}{\nonumber}

\renewcommand{\Im}{\textnormal{Im}\,}

\newcommand{\ZZ}{\mathbb{Z}}
\newcommand{\toa}{\twoheadrightarrow}

\begin{document}
\centerline{\bf\Large Bernoulli free-boundary
problems}\centerline{\bf\Large in strip-like domains}\centerline{\bf
\Large and a property of permanent waves}\centerline{\bf\Large in
water of finite depth } \centerline{}\centerline{\large Eugen
Varvaruca}\centerline{} \centerline{Department of Mathematical
Sciences, University of Bath} \centerline{Claverton Down, Bath BA2
7AY, United Kingdom}\centerline{Email address: {\tt
mapev@maths.bath.ac.uk}}

\begin{abstract} We study weak solutions for a class of free boundary problems which
includes as a special case the classical problem of traveling waves
on water of finite depth. We show that such problems are equivalent
to problems in fixed domains and study the regularity of their
solutions. We also prove that in very general situations the free
boundary is necessarily the graph of a function.
\end{abstract}

\section{Introduction}

One of the classical problems of nonlinear hydrodynamics is that of
traveling two-dimensional gravity waves on water of finite depth,
which arises from the following physical situation. A wave of
permanent form moves with constant speed on the surface of an
incompressible irrotational flow, the bottom of the flow domain
being horizontal. With respect to a frame of reference moving with
the speed of the wave profile, the flow is steady and occupies a
fixed region $\Om$ in the $(X,Y)$-plane, which lies between the real
axis $\mcb$ and some a priori unknown free surface
$\mathcal{S}:=\{(u(s),v(s)):s\in\mathbb{R}\}$. Since the fluid is
incompressible and irrotational, the flow can be described by a
stream function $\psi$ which is harmonic in $\Om$ and satisfies the
following boundary conditions:
\begin{subequations}\label{pb}
\begin{align}
& \psi= constant\quad \text{on }\mcb,\label{p1}\\&\psi=constant\quad
\text{on }{\mathcal S},
\label{p2}\\&\vert\nabla\psi\vert^{2}+2gY=constant \quad\text{on }
{\mathcal S}\label{p3}.
 \end{align}
\end{subequations}
The problem consists of determining the curves $\mcs$ for which a
function $\psi$ with these properties exists in $\Om$.

 There are two particular types of waves
which have received considerable attention in the literature:
periodic waves, where $\mcs$ is assumed to be periodic in the
horizontal direction, and solitary waves, where $\mcs$ is assumed to
be asymptotic to a horizontal line at infinity. Nowadays, the
mathematical theory of these problems contains a wealth of results,
of which some of the most notable are the existence results of Amick
and Toland \cite{AT1, AT2} for both periodic and solitary waves, and
the results on symmetry and monotonicity of solitary waves of Craig
and Sternberg \cite{CS}. We refer the reader to these papers and
their references for historical background and further results. The
related problem of traveling waves with vorticity  \cite{CoS} will
not be addressed in this article.

It is always assumed in the literature that the wave profile $\mcs$
is a graph giving the vertical coordinate as a function of the
horizontal coordinate. But it is more natural to assume a priori
only that $\mcs$ is a curve in parametric form, and to ask whether
there can exist situations in which $\mcs$ is not such a graph. This
is one of the main questions we address in this article, showing
that in very general situations $\mcs$ is necessarily a graph. The
absence up to now of any investigation of this problem is
surprising, since the related problem for waves of infinite depth
was solved long time ago by Spielvogel \cite{Spiel}, in the periodic
setting. He derived a differential inequality for a function $\te$
which gives the angle between the free boundary and the horizontal,
and suggested a geometric argument to show that this  inequality
would prevent the curve from overturning. It is nowadays understood
that for a smooth periodic wave the non-overturning property follows
easily from this inequality combined with the periodicity of $\te$,
see Toland \cite{ToPseudo}. Such an approach can easily be adapted
to the case of smooth periodic waves on water of finite depth, but
it does not work for solitary waves, nor for the wave profiles of
main interest in this article, which may have more general
geometries and need not be smooth.

It is customary in water-wave theory to assume that $\mcs$ is smooth
enough, for example that $\mcs$ is a $C^1$ curve. However, it is
well known that there exist situations in which this smoothness
requirement is not satisfied, namely for `waves of extreme form',
see \cite{AFT, P82}. In recent years, there has been some interest
in weak solutions of the problem of periodic waves of infinite
depth, through a series of papers of Shargorodsky and Toland,
culminating with the comprehensive treatment in \cite{ST}. They
noted that the classical theory of Hardy spaces \cite{Du, Ko} can be
used to assign boundary values to harmonic functions in domains
whose boundary is a locally rectifiable curve, and required the
boundary condition (\ref{p3}) to be satisfied in such a weak sense.
In fact, \cite{ST} deals with and provides a rich mathematical
theory for a more general class of problems, called Bernoulli
problems. Further aspects of this theory, including geometric
properties of free boundaries and the nature of their singularities,
were examined in \cite{EV,EV2}.

In this article we propose an analogue of the theory in \cite{ST}
for a class of problems which would naturally generalize the waves
of finite depth problem. Basically, we keep the boundary conditions
(\ref{p1})-(\ref{p2}), replace (\ref{p3}) by the more general
condition \be \vert\nabla\psi\vert= h(Y)\quad \text{ on }\mcs,\ee
where $h$ is a given function, and provide a interpretation in a
weak sense of these boundary conditions. We impose minimal
smoothness requirements on the curve $\mcs$ and the harmonic
function $\psi$, while the assumptions on the geometry of $\mcs$ are
sufficiently general that, when specialized to the case of water
waves, include both the case of periodic waves and that of solitary
waves.

In \cite{ST} it is assumed that $h$  is a continuous function with
values in $[0,\infty]$ (and avoiding at least one of the values $0$
and $\infty$), which is suitably smooth on the open set where it is
non-zero and finite. In that situation, the singularities of the
curve $\mcs$ can only occur at \emph{stagnation points}, which are
points $(X,Y)$ on $\mcs$ for which $h(Y)=0$ or $h(Y)=\infty$. The
set of stagnation points, denoted by $\msn$, is closed and has
measure zero on $\mcs$, see \cite{ST}, and $\mcs\setminus\msn$ is a
union of smooth open arcs. Aiming for a theory of greater
generality, in this article we consider Bernoulli problems for
functions $h$ which are only Borel measurable. (We believe that it
should be possible to construct examples of explicit solutions of
Bernoulli problems for functions $h$ which are nowhere continuous,
and where $\mcs$ is locally rectifiable but no arc of $\mcs$ is of
class $C^1$.) We note however that the regularity results in
\cite{EV} extend to the present setting, showing that the smoothness
of $h$ on some open intervals on which it is non-zero and finite
implies the regularity of certain arcs of $\mcs\setminus\msn$.

Since one of the main difficulties of Bernoulli problems is the fact
that the free boundary $\mcs$ is unknown a priori, it is of interest
to see if they can be reduced to problems in fixed domains. This is
usually achieved in the hydrodynamics problem by conformal mapping,
but in the situation considered here when $\mcs$ is a curve in
parametric form (and not necessarily a graph) this is not an
entirely straightforward matter, and we are not aware of any
rigorous proof in the literature. We show here that any Bernoulli
free boundary problem is indeed equivalent by means of a conformal
mapping to a problem in a fixed domain, namely a strip.

The main result of the article is that, when $h$ is strictly
decreasing and $\log h$ is concave on the interval where $h$ is
non-zero, any free boundary $\mcs$, possibly with many
singularities, is necessarily the graph of a function. This result
applies in particular to the water-wave problem. In the proof, we
first derive a differential inequality generalizing that of
Spielvogel \cite{Spiel}, and for this purpose the specific Hardy
spaces in which the original problem was posed play an essential
role. To get to the required conclusion, one then uses  a geometric
argument similar in spirit to that first used in \cite{EV2} for
Bernoulli problems in the setting of \cite{ST}. The details here are
however significantly different from those in \cite{EV2}, since the
a priori admissible geometries of the curve $\mcs$ are much more
general here.

 One may wonder whether
any assumptions at all on $h$ are necessary to ensure that $\mcs$ is
a graph. In the setting considered in \cite{ST}, examples of free
boundaries which are not graphs have been constructed in \cite{EV3},
and it seems highly likely that such examples exist in the setting
considered here too.

\section{Bernoulli Free-Boundary Problems}

\subsection{Some Preliminaries}

In this article, by a strip-like domain it is meant a domain
$\Omega$ whose boundary consists of the real axis $\mcb$ and a
non-self-intersecting  curve
$\mathcal{S}:=\{(u(s),v(s)):s\in\mathbb{R}\}$  contained in the
upper half-plane, such that \label{pbx}
 \bese \begin{align}
&\mcs \text{ is locally rectifiable, }\\
 \lim_{s\to\pm\infty} u(s)&=\pm\infty,\quad v\text{ is bounded above}.\end{align}\ese
Let $h:J\to [0,\infty]$ be a given Borel measurable function, where
$J\subset (0,\infty)$ is an interval.
  A Bernoulli problem is one of finding a strip-like domain
$\Omega$ in which there exists a harmonic function $\psi$ satisfying
the following boundary conditions:
\begin{align}
\psi&=0 \quad\quad\text{ almost everywhere on } \mcb\label{pba},\\
\psi&=C\quad\quad\text{ almost everywhere on } \mcs,\label{pbb}
\\\vert\nabla\psi\vert&= h(Y)\,\,\,\quad \text{almost everywhere on }
{\mathcal S}\label{pbc},
 \end{align}
where $C$ is a positive constant. (Whenever we refer to a relation
which holds on a locally rectifiable curve, such as in
(\ref{pba})-(\ref{pbc}), `almost everywhere' refers to
one-dimensional Hausdorff measure, or arclength.) Since the curve
${\mathcal S}$ is not prescribed a priori, it is called a \emph{free
boundary}.

We now explain the weak sense in which these conditions are to be
satisfied. In short, if the following conditions hold: \be
\psi\text{ is bounded in }\Om,\label{psb}\ee \be \text{the
subharmonic function $|\nabla\psi|$ has a harmonic majorant in
$\Om$,}\label{hamj}\ee then the function $\psi$ and its partial
derivatives have non-tangential limits almost everywhere on $\mcs$
and $\mcb$. In this situation, the required conditions
(\ref{pba})-(\ref{pbc}) would refer to the non-tangential boundary
values of $\psi$ and $\nabla\psi$.

To support the above claim, we need a precise definition of a
non-tangential limit, and the definition and properties of Hardy
spaces in general domains.

\subsection{Non-tangential Limits}

Let $\Th$ be a bounded open set in the plane whose boundary is a
rectifiable Jordan curve $\mathcal{T}$. Let $(X_0, Y_0)$ be a point
on $\mct$ at which $\mct$ has a tangent, and let $\mathbf{n}_i$ be
the corresponding unit inner normal vector. We say that a sequence
$\{(X_n, Y_n)\}\indn$ of points in $\Th$ \emph{tends to $(X_0, Y_0)$
non-tangentially within $\Th$} if $(X_n, Y_n)\to (X_0,Y_0)$ as
$n\to\infty$ and there exists $\kappa>0$ such that
\[(X_n-X_0, Y_n-Y_0)\cdot\mathbf{n}_i\geq \kappa[(X_n-X_0)^2+(Y_n-Y_0)^2]^{1/2}
\quad\text{for all }n\geq 1,\] where the dot denotes the usual inner
product in $\bdr^2$.

 Let $F:\Th\to\bdc$ be a function,
and let $l\in\bdc$. We say that $F$ \emph{has non-tangential limit
$l$ at $(X_0,Y_0)$ within $\Th$}, and write
\[\lim_{(X,Y)\toa(X_0, Y_0)} F(X,Y)=l,\]
if $\lim_{n\to\infty}F(X_n, Y_n)=l$ for every sequence $\{(X_n,
Y_n)\}\indn$  which tends to $(X_0, Y_0)$ non-tangentially within
$\Th$.

Let $\Om$ be a strip-like domain, let $Z_0:=X_0+iY_0$ be a point on
$\mcs$ at which $\mcs$ has a tangent, and let $F:\Om\to\bdc$. Let
$Z_1$ and $Z_2$ be two points on $\mcs$ such that $Z_0$ is between
$Z_1$ and $Z_2$. Then $Z_1$ and $Z_2$ can be joined by a rectifiable
arc contained in $\Om$. This arc, together with the arc $\mca$ of
$\mcs$ joining $Z_1$ and $Z_2$, determines a rectifiable Jordan
curve, which is the boundary of a bounded subdomain $\Th$ of $\Om$.
Let $l\in\bdc$. We say that $F$ \emph{has non-tangential limit $l$
at $(X_0,Y_0)$ within $\Om$} if $F$ has non-tangential limit $l$ at
$(X_0,Y_0)$ within $\Th$. It is easy to see that this definition is
meaningful, in the sense that it does not depend on the set $\Th$
used.

\subsection{Hardy Spaces}

Let $D$ be the unit disc in the complex plane. For $p\in[1,\infty)$,
the Hardy space $h^p_\bdc(D)$ is usually defined, see \cite{Du,Ko,
Ru}, as the class of harmonic functions $F:D\to\bdc$  with the
property that \be\sup_{r\in
(0,1)}\int_{-\pi}^\pi|F(re^{it})|^p\,dt<+\infty.\label{hagr}\ee The
Hardy space $h^\infty_\bdc(D)$ is the class of bounded harmonic
functions in $D$. For $p\in[1,\infty]$, the Hardy space $\ha^p(D)$
is the class of holomorphic functions in $\h^p(D)$.  It is well
known that any function in $h_\bdc^p(D)$, $p\in [1,\infty]$, has
non-tangential limits almost everywhere on the unit circle. The M.
Riesz Theorem \cite[Theorem 4.1, p.\ 54]{Du} asserts that, if $U\in
\h^p(D)$ for some $p\in (1,\infty)$, and if $V$ is a harmonic
function such that $U+iV$ is holomorphic, then $V\in \h^p(D)$.

 An important fact, which leads to the definition of Hardy spaces in general
domains \cite[Ch.\ 10]{Du}, is that, for $p\in [1, \infty)$, a
harmonic function $F$ belongs to $h_\bdc^p(D)$ if and only if the
subharmonic function $|F|^p$ has a harmonic majorant, i.e. there
exists a positive harmonic function $u$ in $D$ such that $|F|^p\leq
u$ in $D$.  Let $\Th$ be an open set. For $p\in [1,\infty)$, the
space $\h^p(\Th)$ is  the class of harmonic functions $F:\Th\to\bdc$
for which the subharmonic function $|F|^p$ has a harmonic majorant
in $\Th$. The Hardy space $\h^\infty(\Th)$ is the class of bounded
harmonic functions in $\Th$. The spaces $\ha^p(\Th)$ consists of the
holomorphic functions in $\h^p(\Th)$, for $p\in[1,\infty]$. It is
easy to check that the Hardy spaces are conformally invariant: if
$\Th_1$ and $\Th_2$ are two open sets, and $\gamma:\Th_1\to\Th_2$ is
a conformal mapping, then $F\in\h^p(\Th_2)$ if and only if
$F\circ\gamma\in \h^p(\Th_1)$, where $p\in [1,\infty]$. Due to this
fact, many properties of the Hardy spaces of the disc extend by
conformal mapping to Hardy spaces of simply connected domains.
 If $\Th$ is a bounded domain whose
boundary is a rectifiable Jordan curve, then any function in
$\h^p(\Th)$, where $1\leq p\leq\infty$, has non-tangential boundary
values almost everywhere. It is immediate that, if $\Om$ is a
strip-like domain, then any function in $\h^p(\Om)$, where $1\leq
p\leq\infty$, has non-tangential boundary values almost everywhere
on $\mcs$ and $\mcb$.

Finally, we mention that it is possible to define Hardy spaces
$\ha^p(\Th)$ also for values $p\in (0,1)$, namely as the class of
holomorphic functions $F:\Th\to\bdc$ for which the subharmonic
function $|F|^p$ has a harmonic majorant in $\Th$, an open set.
(Note that when $F$ is only harmonic in $\Th$ and $p\in(0,1)$, the
function $|F|^p$ need not be subharmonic.) For the unit disc $D$,
the class $\ha^p(D)$, $p\in(0,1)$, coincides with the class of
holomorphic functions for which (\ref{hagr}) holds. In this article,
only marginal use is made of the spaces $\ha^p(\Th)$, $p\in(0,1)$.

\subsection{Further Preliminaries}

If $\psi$ is harmonic in a strip-like domain $\Om$ and satisfies
(\ref{psb}) and (\ref{hamj}), the preceding discussion ensures that
$\psi$ and its partial derivatives have non-tangential boundary
values almost everywhere on $\mcs$ and $\mcb$. In a Bernoulli
problem, these boundary values are required to satisfy
(\ref{pba})-(\ref{pbc}).

 Whenever $\psi$ is a harmonic function in a strip-like domain $\Om$, we shall consider also a harmonic function
$\vfi$ such that the function $\om:=\vfi+i\psi$ is holomorphic in
$\Om$. The condition (\ref{hamj}) can then be reformulated as
$\om'\in\ha^1(\Om)$. Here and in what follows, $'$ denotes
differentiation (it will be clear from the context whether either
real or complex differentiation is meant).
  Let also $\Om^R$ denote the reflection
of $\Om$ in the line $\mcb$, let \be
\widetilde\Om:=\Om\cup\mcb\cup\Om^R,\label{tio}\ee
 and let $L>0$ be such that $\widetilde\Om$ is contained in the strip
determined by the lines $Y=-L$ and $Y=L$.

The next result gathers several properties of weak solutions of
Bernoulli problems, refining some results in \cite{ST}. Whilst in
\cite{ST} the function $\psi$ was assumed continuous in
$\Om\cup\mcs$, here we derive the continuity of $\om$ in
$\Om\cup\mcs$ from (\ref{psb}) and (\ref{pbb}). Also, the claims in
\cite{ST} that the condition (\ref{pbb}) ensures that formally the
tangential derivative of $\psi$ on $\mcs$ is zero, and hence that
the condition (\ref{pbc}) is equivalent to a Neumann condition, are
put on a rigorous basis here. (The fact proved here that $\om$ has a
holomorphic extension to $\widetilde\Om$ has no analogue in the
setting of \cite{ST}.)

\begin{proposition}\label{PPR} Let $\Om$ be a strip-like domain, and let $\om=\vfi+i\psi$
 be a holomorphic function
in $\Om$ such that {\rm(\ref{psb})} holds.
\begin{itemize}
\item[(i)] If  $\psi$ satisfies
{\rm(\ref{pba})}, then $\om$ has a holomorphic extension to
$\widetilde\Om$ which satisfies $\om(\bar Z)=\overline{\om(Z)}$ for
all $Z\in\widetilde\Om$.
\item[(ii)] If $\psi$ satisfies {\rm(\ref{pbb})}, then $\om$ has a
continuous extension to $\Om\cup\mcs$.

\item[(iii)] If {\rm(\ref{hamj})} holds and $\psi$ satisfies {\rm(\ref{pbb})}, then
\be\Big(\lim_{(X,Y)\toa
(X_0,Y_0)}\nabla\psi(X,Y)\Big)\cdot\mathbf{t}(X_0,Y_0)=0\quad
\text{a.e.\ $(X_0, Y_0)$ on }\mcs,\label{tanz}\ee where
$\mathbf{t}(X_0,Y_0)$ is a unit tangent vector to $\mcs$ at
$(X_0,Y_0)$.

\item[(iv)] If {\rm(\ref{hamj})} holds and $\psi$ satisfies {\rm(\ref{pbb})},
 then $\psi$ satisfies {\rm(\ref{pbc})} if and only if
\be \Big(\lim_{(X,Y)\toa
(X_0,Y_0)}\nabla\psi(X,Y)\Big)\cdot\mathbf{n}_o(X_0,Y_0)=h(Y)\quad\text{a.e.\
$(X_0, Y_0)$ on }{\mathcal S}\label{pbij},\ee where
$\mathbf{n}_o(X_0,Y_0)$ is the unit outer normal to $\Om$ at
$(X_0,Y_0)$.
\end{itemize}
\end{proposition}

It is necessary for our purposes to strengthen the basic conditions
(\ref{psb})-(\ref{hamj}). The formal definition of a Bernoulli
problem given below involves a holomorphic function $\om=\vfi+i\psi$
in the domain $\widetilde\Om$ given by (\ref{tio}). The conditions
(\ref{a1})-(\ref{a2}) are motivated by Proposition \ref{PPR}.

\begin{definition} A Bernoulli free-boundary problem, or problem {\bf(B)}, is to find a strip-like
domain $\Om$ and a function $\om= \phi+i\psi$ in $\widetilde\Om$,
such that \be\text{$\om$ is holomorphic in $\widetilde\Om$ and
continuous on its closure},\label{a1} \ee \be \om(\bar
Z)=\overline{\om(Z)}\quad\text{for all
}Z\in\widetilde\Om,\label{a2}\ee \be
\om'\in\ha^1(\widetilde\Om),\label{a3}\ee \be \om'\neq 0\text{ in
}\widetilde\Om\quad\text{and}\quad
1/\om'\in\ha^1(\widetilde\Om),\label{a4}\ee and {\rm (\ref{pbb})},
{\rm (\ref{pbc})} hold.
\end{definition}

\begin{remark}\label{rd}{\rm It will be seen in the proof of Theorem \ref{teqi} that
any solution $\om$ of (\ref{a1}), (\ref{a2}) and (\ref{pbb})
automatically satisfies the condition $\om'\neq 0$ in
$\widetilde\Om$  required in (\ref{a4}).}
\end{remark}

\section{The Main Results}

\subsection{An Equivalent Problem in a Strip}

We now introduce a problem in a strip, to which we shall prove that
problem $\mathbf{(B)}$ is equivalent. Firstly, let us introduce some
more notation. For any $a,b\in\bdr$ with $a<b$, let us denote \[
\Pi_{a,b}:=\{x+iy\in\bdc:a<y<b\},\] and for any $c\in\bdr$, let \[
\mcl_c:=\{x+iy\in\bdc:y=c\}.\]

\begin{definition}
We say that a function $W=U+iV$ in the strip $\Pi_{-C,C}$ is a
solution of problem {\bf(P)} if the following conditions are
satisfied: \be\text{$W$ is holomorphic in  $\Pi_{-C,C}$ and
continuous on its closure},\label{holw}\ee \be W(\bar
z)=\overline{W(z)}\quad\text{for all
}z\in\Pi_{-C,C},\label{c4}\ee\be V \text{ is bounded},\qquad V(t,
C)>0\quad\text{for all }t\in\bdr\label{c1},\ee \be \text{the mapping
}t\mapsto W(t,C)\text{ is injective and
}\lim_{t\to\pm\infty}U(t,C)=\pm\infty, \label{c3}\ee \be W'\in
\ha^1(\Pi_{-C,C}),\label{c2} \ee \be W'\neq 0\text{ in }\Pi_{-C,
C}\quad\text{and}\quad 1/W'\in \ha^1(\Pi_{-C,C}),\label{c7}\ee
   \be h(V(t,C))|\nabla
V(t,C)|=1\quad\text{for almost every }t\in\bdr.\label{c6}\ee

\end{definition}

\begin{remark}
{\rm The condition (\ref{c2}) ensure that the partial derivatives of
$V$ have non-tangential boundary values almost everywhere. The
condition (\ref{c6}) refers to these non-tangential boundary values
on $\mcl_C$.}
\end{remark}

\begin{remark}\label{rtt} {\rm It will be seen in the proof of Theorem
\ref{teqi} that any solution $W$ of {\rm
(\ref{holw})}-{\rm(\ref{c3})} necessarily satisfies  the condition
that $W'\neq 0$ in $\Pi_{-C,C}$ required in (\ref{c7}).}
\end{remark}

The main result on equivalence is the following.

\begin{theorem}\label{teqi} Let $(\Om,\om)$ be a solution of problem {\bf
(B)}. Then $\om$ is a conformal mapping from $\widetilde\Om$ onto
the strip $\Pi_{-C, C}$ and a homeomorphism from the closure of
$\widetilde\Om$ onto the closure of $\Pi_{-C,C}$, with $\Om$ being
mapped onto $\Pi_{0,C}$, $\mcb$ being mapped onto $\mcl_0$ and
$\mcs$ being mapped onto $\mcl_{C}$. Let $W$ be the inverse
conformal mapping, from $\Pi_{-C,C}$ onto $\widetilde\Om$. Then $W$
is a solution of problem {\bf (P)}.

Conversely, let $W$ be a solution of problem {\bf (P)}.  Let
$\mcs:=\{W(t,C):t\in\bdr\}$, which is a non-self-intersecting curve,
and let $\Om$ be the domain whose boundary consists of the curve
$\mcs$ and the real axis $\mcb$. Then $\Om$ is a strip-like domain.
Let $\widetilde\Om$ be given by {\rm (\ref{tio})}. Then $W$ is a
conformal mapping from $\Pi_{-C,C}$ onto $\widetilde\Om$ and a
homeomorphism from the closure of $\Pi_{-C,C}$ onto the closure of
$\widetilde\Om$, with  $\Pi_{0,C}$ being mapped onto $\Om$, $\mcl_0$
being mapped onto $\mcb$ and $\mcl_C$ being mapped onto $\mcs$. Let
$\om$ be the inverse conformal mapping, from $\widetilde\Om$ onto
$\Pi_{-C,C}$. Then $(\Om,\om)$ is a solution of problem {\bf (B)}.

\end{theorem}

\subsection{Local Regularity}

The following local regularity result is an analogue of
\cite[Theorem 2.1 and Theorem 2.3]{EV}.

\begin{theorem}\label{tlre}
Let $h:(0,\infty)\to [0,\infty]$ be a Borel measurable function. Let
$(c,d)\subset (0,\infty)$ be an interval, and suppose that
\begin{align} h&\in C^{n,\alpha}_{\textnormal
{loc}}((c,d))\text{ where }n\in \mathbb{N}\cup\{0\},\,\, \alpha\in
(0,1),\non\\
h&\neq 0\text{ on }(c,d).\non
 \end{align}
Let $(\Om,\om)$ be a solution of problem {\bf (B)}. Let $Z_1,
Z_2\in\mcs$ with $Z_1\neq Z_2$ be such that the open arc
$\mathcal{A}$ of $\mcs$ joining $Z_1$ and $Z_2$ is contained in the
strip $\Pi_{c,d}$.  Then $\mca$ is a curve of class
$C^{n+1,\alpha}_{\textnormal {loc}}$ and $\om\in
C^{n+1,\alpha}_{\textnormal {loc}}(\Om\cup\mathcal{A})$.

\end{theorem}

\subsection{When the Free Boundary is a Graph}

The main result of this article is the following.

\begin{theorem}\label{tgraph}
Let $d>0$ and $h:(0,d]\to[0,\infty)$ be such that \bese\begin{align}
&h(d)=0\text{ and } h>0\text{ on }(0,d),\\ &\text{$h\in C((0,d])\cap
C^{1,\alpha}_{\textnormal{loc}}((0,d))$
for some $\alpha\in (0,1)$},\\
 &\text{$h$ is strictly decreasing and $\log h$ is concave on $(0,d)$}.
\label{nash}\end{align}\ese Let $(\Om,\om)$ be a solution of problem
{\bf (B)}. Then there exists a parametrization $\{(u(t), v(t)):
t\in\bdr\}$ of $\mcs$,  where $u$, $v$ are locally absolutely
continuous and \be v'(t)>0\quad\text{for almost all }t\in\bdr.
\label{vepo}\ee Hence there exists a continuous function $\eta:\bdr
\to\bdr$ such that \[ \mathcal{S}=\{(X,\eta(X)):\,\,X\in
 \mathbb{R}\}. \]

\end{theorem}

\section{Proofs}

In this section we give the proofs of the results stated in the
previous sections. We start with a lemma which will be useful for
some of the proofs. The first part is well known, while the second
part is essentially a local version of Privalov's Theorem, see
\cite[Lemma 2.2]{EV} for the proof of a very similar result.

\begin{lemma}\label{lreg} Let $x_0\in\bdr$  and $r>0$. For $t>0$,
let us denote
$B^+_t(x_0):=\{(x,y)\in\bdr^2:(x-x_0)^2+y^2<t^2,\,y>0\}$. Let
$F\in\ha^1(B_r^+(x_0))$ be of the form $F=U+iV$, where $U$ and $V$
are real-valued functions, and let, for $s\in (x_0-r,x_0+r)$,
\[u(s)+iv(s):=\lim_{(x,y)\toa (s,0)}\big(U(x,y)+iV(x,y)\big).\]
\begin{itemize}
\item[(i)] If $u$ is continuous on $(x_0-r,x_0+r)$, then $U$ is
continuous on $B^+_r(x_0)\cup \{(x,0):x\in(x_0-r,x_0+r)\}$.
\item[(ii)] If $u\in C^{k,\al}_{\textnormal{loc}}((-r,r))$ for some
$k\in\mathbb{N}\cup\{0\}$ and $\alpha\in (0,1)$, then $v\in
C^{k,\al}_{\textnormal{loc}}((-r,r))$ and $U$, $V$ are of class
$C^{k,\alpha}$ in the closure of $B_t^+(x_0)$ for every $t\in
(0,r)$.
\end{itemize}

\end{lemma}
\subsection{Proof of Proposition \ref{PPR}}

\begin{proof}[Proof of Proposition \ref{PPR}]
(i) Since $\psi$ is bounded in $\Om$, it is immediate from the first
part of Lemma \ref{lreg} that $\psi$ extends continuously to
$\Om\cup\mcb$. The conclusion now follows from the classical Schwarz
Reflection Principle, see \cite[Theorem 11.14, p.\ 237]{Ru}.

(ii) Let $Z_1$ and $Z_2$ be any two distinct points on $\mcs$. Then
they can be joined by a rectifiable arc contained in $\Om$. This
arc, together with the arc $\mathcal{A}$ of $\mcs$ joining $Z_1$ and
$Z_2$, determines a rectifiable Jordan curve, which is the boundary
of a bounded subdomain $\Th$ of $\Om$. Let $\gamma:D\to\Th$ be a
conformal mapping from $D$ onto $\Th$. By Caratheodory's Theorem,
see \cite[Ch.\ II, \S C, p.\ 35]{Ko}, $\gamma$ is a homeomorphism
from the closure of $D$ onto the closure of $\Th$. Let
$I\subset\bdr$ be an interval such that
$\gamma(\Gamma_I)=\mathcal{A}$, where $\Gamma_I:=\{e^{is}:s\in I\}$.

Let $F:=\om\circ\gamma$. Then $F\in\ha^1(D)$. Moreover, $F-C$ is
real-valued on $\Gamma_I$, so by the Schwarz Reflection Principle
for $\ha^1(D)$ functions, see \cite[Ch.\ III, \S E, p.\ 63]{Ko}, $F$
has a holomorphic extension across $\Gamma_I$. In particular $F$
extends continuously to $D\cup\Gamma_I$, and therefore $\om$ has a
continuous extension to $\Om\cup\mathcal{A}$. Since $Z_1$ and $Z_2$
were arbitrary on $\mcs$, this proves (ii).

(iii) Suppose now that (\ref{hamj}) holds. Let $\Th$ be as in the
proof of (ii). Since the boundary of $\Th$ is a rectifiable Jordan
curve, it follows from \cite[Theorem 3.11 and Theorem 3.12, p.\
42]{Du} that $\gamma'\in\ha^1(D)$, the mapping $s\mapsto
\gamma(e^{is})$ is absolutely continuous and
\be\frac{d}{ds}\gamma(e^{is})=ie^{is}\lim_{\zeta\toa e^{is}}\gamma
'(\zeta)\quad\text{almost everywhere}. \label{der}\ee Note that the
mapping $I\ni s\mapsto \gamma(e^{is})$ is an absolutely continuous
parametrization of $\mathcal{A}$ and, almost everywhere on $I$,
$\frac{d}{ds}\gamma(e^{is})$ gives a tangent vector to $\mcs$.

 Since $F$ is holomorphic across $\Gamma_I$ and $F-C$ is real-valued on
 $\Gamma_I$, it follows that
\[\Im \{ie^{is}F'(e^{is})\}=0\quad\text{for all }s\in I.\]
Since $F=\om\circ\gamma$ on $D$, the above means that
\[\Im\{ie^{is}\lim_{\zeta\toa e^{is}}(\omega'(\gamma(\zeta))\gamma'(\zeta))\}=0\quad\text{ almost everywhere on }I.\]
It follows, using (\ref{der}), that
\[\Im\left\{\Big(\lim_{X+iY\toa \gamma(e^{is})} \omega'(X+iY)\Big)\frac{d}{ds}\gamma(e^{is})\right\}=0\quad\text{ almost everywhere on }I,\]
and therefore
\[\Big(\lim_{(X,Y)\toa (X_0, Y_0)}\nabla\psi(X,Y)\Big)\cdot\mathbf{t}(X_0, Y_0)=0 \quad\text{a.e.\ $(X_0,Y_0)$ on }\mathcal{A}.\]
Since $Z_1$ and $Z_2$ were arbitrary on $\mcs$, (\ref{tanz})
follows.

(iv) The Maximum Principle shows that \be
0<\psi<C\quad\text{everywhere in }\Om.\label{inq}\ee Since
(\ref{inq}) holds, it follows that
\be\Big(\lim_{(X,Y)\toa(X_0,Y_0)}\nabla\psi(X,Y)\Big)\cdot\mathbf{n}_o(X_0,Y_0)\geq
0 \quad\text{a.e. $(X_0,Y_0)$ on }\mcs.\label{norp}\ee The required
result easily follows from (\ref{tanz}) and (\ref{norp}).

\end{proof}

\subsection{Hardy Spaces of the Strip}

We now summarize some facts concerning the Hardy spaces of a strip,
a thorough treatment of which has recently been given in \cite{BK}.
However, most (though not all) of the results below follow
immediately by means of conformal mapping from the corresponding
results in the unit disc, which are well known.

 Let $q:\bdr\to\bdr$ be given by $q(t)=e^{-\frac{\pi}{2C}|t|}$ for all
$t\in \bdr$. For any $p\in [1,\infty)$, let \be
L^p(\bdr,q(t)\,dt):=\left\{F:\bdr\to\bdc:\int_{-\infty}^{\infty}|F(t)|^p
q(t)\,dt<+\infty\right\}.\ee If $F$ is in the Hardy space
$\h^p(\Pi_{-C, C})$, where $p\in[1,\infty]$,  then $F$ has
non-tangential boundary values almost everywhere, which, when
$p\in[1,\infty)$,
 satisfy
$|F(\cdot,\pm C)|\in L^p(\bdr,q(t)\,dt)$. For $p\in (1,\infty]$, but
not necessarily when $p=1$, $F$ can be expressed as a Poisson
Integral of its boundary values. Here, for every two functions
 $F_C,F_{-C}\in L^1(q(t)\,dt)$, the Poisson
 Integral \cite{Wi} associated to $F_C$ and $F_{-C}$ is given, for all $(x,y)\in\Pi_{-C,C}$, by
\begin{align}
\mathcal{P}[F_C,F_{-C}](x,y):&=\frac{1}{4C}\int_{-\infty}^{\infty}\frac{\cos(\frac{\pi}{2C}y)}
{\cosh(\frac{\pi}{2C}(x-s))-\sin(\frac{\pi}{2C}y)}\,F_C(s)\,ds\non\\
&+\frac{1}{4C}\int_{-\infty}^{\infty}\frac{\cos(\frac{\pi}{2C}y)}
{\cosh(\frac{\pi}{2C}(x-s))+\sin(\frac{\pi}{2C}y)}\,F_{-C}(s)\,ds.\non
\end{align}
However, if $F\in \ha^1(\Pi_{-C,C})$ then, by the F.\ and M.\ Riesz
Theorem \cite[Ch.\ 2]{Ko}, $F$ can be expressed as the Poisson
Integral of its boundary values.

For any $F\in L^1(\bdr,q(t)\,dt)$, we denote by $\mathcal{P}_e[F]$
and $P_o[F]$ the functions in $\Pi_{-C,C}$ given by \be
\mathcal{P}_e[F]:=\mathcal{P}[F,F]\quad\text{and}\quad
\mathcal{P}_o[F]:=\mathcal{P}[F,-F]. \label{Pois}\ee

The next result is slightly more general than needed for our
purposes, but we think it might be of interest in itself. It is an
analogue of \cite[Theorem 3.11, p.\ 42]{Du}.

\begin{lemma}\label{hac} Let the holomorphic function
$F:\Pi_{-C,C}\to\bdc$ be such that $F'\in \ha^1(\Pi_{-C,C})$. Then
$F$ is continuous in the closure of $\Pi_{-C,C}$, the mappings
$\bdr\ni t\mapsto F(t\pm iC)$ are locally absolutely continuous, and
\[\frac{d}{dt}F(t\pm iC)=\lim_{(x,y)\toa (t,\pm C)} F'(x+iy)\quad\text{for a.e.\ }t\in\bdr.\]
\end{lemma}

\begin{proof}[Proof of Lemma \ref{hac}]
 For any $A>0$, let
\[\mathcal{R}_A:=\{(x,y)\in\Pi_{-C,C}:-A<x<A\}.\]
Let $\gamma:D\to\mathcal{R}_A$ be a conformal mapping. Then $\gamma$
is a homeomorphism from the closure of $D$ onto the closure of
$\mathcal{R}_A$ and, since the boundary of $\mathcal{R}_A$ is a
rectifiable Jordan curve, it follows from \cite[Theorem 3.12, p.\
44]{Du} that $\gamma'\in\ha^1(D)$. Let $G=F\circ\gamma$. Then, since
$G'=F'(\gamma)\gamma'$, we deduce that $G'\in\ha^{1/2}(D)$. By a
classical result of Hardy and Littlewood, see \cite[Theorem 5.12,
p.\ 88]{Du}, it follows  that $G\in \ha^1(D)$, and therefore
$F\in\ha^1(\mathcal{R}_A)$. Since this is true for any $A>0$, it
follows that $F$ has non-tangential boundary values almost
everywhere, which we denote $F(\cdot\pm iC)$. Let us also denote
\[f(t\pm iC):=\lim_{(x,y)\toa (t,\pm C)} F'(x+iy) \quad\text{for a.e.\
}t\in\bdr.\]

 Let us focus attention on the behaviour of $F$ near the line $\mcl_C$.
 It is a consequence of \cite[Corollary 2.1]{BK} that,
for every $A>0$, \be ||F'(\cdot+iy)-f(\cdot+iC)||_{L^1(-A,A)}\to
0\quad\text{as }y\nearrow C.\label{baka}\ee Let $x_0\in\bdr$ be such
that there exists $\lim_{y\nearrow C}F(x_0+iy)=:F(x_0+iC)$. For
every $x\in\bdr$ and $y\in (0,C)$, \be
F(x+iy)-F(x_0+iy)=\int_{x_0}^x F'(s+iy)\,ds.\label{zar}\ee We deduce
from (\ref{zar}), upon passing to the limit as $y\nearrow C$ and
taking into account (\ref{baka}), that for every $x$ such that there
exists $\lim_{y\nearrow C}F(x+iy)=:F(x+iC)$, the following holds \[
F(x+iC)-F(x_0+iC)=\int_{x_0}^x f(s+iC)\,ds.\] Hence $\bdr\ni
x\mapsto F(x+iC)$ coincides almost everywhere with a locally
absolutely continuous function. By the first part of Lemma
\ref{lreg}, $F$ is continuous on $\Pi_{-C,C}\cup\mcl_C$. Hence, the
mapping $\bdr\ni x\mapsto F(x+iC)$ is locally absolutely continuous
with
\[\frac{d}{dt}F(t+ iC)= f(t+iC)\quad\text{for a.e.\ }t\in\bdr.\]

A similar argument can be used to deal with the behaviour of $F$
near the line $\mcl_{-C}$. This completes the proof of Lemma
\ref{hac}.
\end{proof}

\subsection{Proof of Theorem \ref{teqi}}

\begin{proof}[Proof of Theorem \ref{teqi}]
We start by proving the result claimed in Remark \ref{rd}. Let $\om$
satisfy (\ref{a1}), (\ref{a2}) and (\ref{pbb}).
 We first
prove the existence of a conformal mapping $W_0$ from $\Pi_{-C,{C}}$
onto $\widetilde\Om$, such that $W_0(\bar z)=\overline{W_0(z)}$ for
all $z\in\Pi_{-C,C}$, and which has an extension as a homeomorphism
from the closure of $\Pi_{-C,{C}}$ to the closure of
$\widetilde\Om$, with $\mcl_0$ being mapped onto $\mcb$ and
$\mcl_{C}$ being mapped onto $\mcs$. Indeed, let $\beta$ be the
conformal mapping from $D$ onto $\Pi_{-C,C}$ given by
\be\beta(\xi)=\frac{2C}{\pi}\log\frac{1+\xi}{1-\xi}\qquad\text{for
all }\xi\in D. \label{be}\ee  Let $\alpha$ be a conformal map from
$D$ onto the strip containing $\widetilde\Om$ which is determined by
the lines $Y=-L$ and $Y=L$, given by
\be\alpha(\zeta)=\frac{2L}{\pi}\log\frac{1+\zeta}{1-\zeta}\qquad\text{for
all }\zeta\in D. \label{al}\ee Let $\Th=\alpha^{-1}(\widetilde\Om)$.
Then $\Th$ is a subdomain of $D$, whose boundary is a Jordan curve
consisting of two arcs which are symmetric about the real axis,
contained in $D$ and joining the points $-1$ and $1$. An immediate
application of Caratheodory's Theorem ensures the existence of a
conformal mapping $\delta$ from $D$ onto $\Th$ which extends as a
homeomorphism between the closures of these domains, and is such
that $\delta(\pm 1)=\pm 1$ and $\delta$ maps the segment $[-1,1]$ of
the real line onto itself. Denoting
$W_0:=\alpha\circ\delta\circ\beta^{-1}$, this mapping $W_0$ has all
the required properties. It follows that $\psi\circ W_0$ is bounded
in $\Pi_{-C,{C}}$ and satisfies $\psi\circ W_0=C$ on $\mcl_{C}$ and
$\psi\circ W_0=-C$ on $\mcl_{-C}$. It follows from the Maximum
Principle that $\psi(W_0(x,y))=y$ in $\Pi_{-C,C}$, and therefore
$\vfi(W_0(x,y))=x+c_0$ in $\Pi_{-C,C}$. Hence $\om-c_0$ is the
inverse of $W_0$, and therefore $\omega'\neq 0$ in $\widetilde\Om$.
This proves the result claimed in Remark \ref{rd}.

Suppose now that $(\Om,\om)$ is a solution of problem {\bf (B)}. The
preceding considerations show that $\om$ is a conformal mapping from
$\widetilde\Om$ onto the strip $\Pi_{-C, C}$ and a homeomorphism
from the closure of $\widetilde\Om$ onto the closure of
$\Pi_{-C,C}$. It is immediate that the inverse conformal mapping $W$
satisfies (\ref{holw})-(\ref{c7}). Also, an argument similar to that
in \cite[Section 3.5, p.\ 43]{Du} shows that, for any Borel set
$\mca$ of $\bdr$, $\mca$ has measure zero if and only if $W(\mca)$
has one-dimensional Hausdorff measure zero on $\mcs$. Moreover,
almost every $t_0\in\bdr$ has the following property: a sequence
$\{(x_n,y_n)\}_{n\geq 1}$ tends to $(t_0,C)$ non-tangentially within
$\Pi_{0,C}$ if and only if the sequence $\{W(x_n,y_n)\}_{n\geq 1}$
tends to $W(t_0,C)$ non-tangentially within $\Om$. In view of these
facts, (\ref{c6}) follows from (\ref{pbc}).  This completes the
proof of the fact that $W$ is a solution of problem {\bf (P)}.

We now prove the result claimed in Remark \ref{rtt}. Let $W$ satisfy
(\ref{holw})-(\ref{c3}). Let $\mcs:=\{W(t,C):t\in\bdr\}$, which by
the first part of (\ref{c3}) is a non-self-intersecting curve. Let
$\Om$ be the domain whose boundary consists of the curve $\mcs$ and
the real axis $\mcb$, and let $\widetilde\Om$ be given by {\rm
(\ref{tio})}. The first part of (\ref{c1}) ensures, by means of the
M. Riesz Theorem, that $U$ can be recovered as the Poisson Integral
of its boundary values in the strip $\Pi_{-C,C}$. We deduce from
this that \be \lim_{x\to\pm\infty} U(x,y)=\pm\infty
\quad\text{uniformly in }y\in[-C,C].\label{liu}\ee Consider again
the mappings $\beta$ and $\alpha$ given by (\ref{be}) and
(\ref{al}). Let $\sigma:=\alpha^{-1}\circ W\circ\beta$. Then
$\sigma$ is holomorphic in $D$, and the continuity of $W$ in the
closure of $\Pi_{-C,C}$ and its behavior expressed by (\ref{liu})
ensure that $\sigma$ has a continuous extension to the closure of
$D$. Moreover, (\ref{c3}), (\ref{c4}) and (\ref{c1})  show that
$\sigma$ is injective on the boundary of $D$. It follows from the
 classical Darboux-Picard Theorem \cite[Corollary 9.16, p.\ 310]{Bu}
 that $\sigma$ is a conformal bijection from
$D$ onto $\Th$, where $\Th=\alpha^{-1}(\widetilde\Om)$, a subdomain
of $D$ whose boundary is a Jordan curve consisting of two arcs which
are symmetric about the real axis, contained in $D$ and joining the
points $-1$ and $1$. It now follows that $W$ is a conformal mapping
from $\Pi_{-C,C}$ onto $\widetilde\Om$, and therefore $W'\neq 0$ on
$\Pi_{-C, C}$, as required.

Suppose now that $W$ is a solution of problem {\bf (P)}. It follows
from (\ref{c2}) upon invoking Lemma \ref{hac} that the mapping
$\bdr\ni t\mapsto W(t,C)$ is locally absolutely continuous, and
hence $\mcs$ is a locally rectifiable curve. Therefore $\Om$ is a
strip-like domain. We have already seen that $W$ is a conformal
mapping from $\Pi_{-C,C}$ onto $\widetilde\Om$ and a homeomorphism
from the closure of $\Pi_{-C,C}$ onto the closure of
$\widetilde\Om$. It is immediate that, if $\om$ is the inverse
conformal mapping, then (\ref{a1})-(\ref{a4}) and (\ref{pbb}) hold.
The same argument used at the end of the first part of the proof
shows that (\ref{pbc}) follows from (\ref{c6}). This completes the
proof of the fact that $(\Om,\om)$ is a solution of problem {\bf
(B)}.
\end{proof}

\subsection{Sketch of the Proof of Theorem \ref{tlre}}

\begin{proof}[Proof of Theorem \ref{tlre}]
The proof is based on arguments which are very similar to some used
in \cite{EV}. Because of Theorem \ref{teqi}, one can concentrate on
proving the regularity of the corresponding solution of problem~{\bf
(P)}.

We start by deriving some further properties of solutions $W$ of
problem {\bf (P)}. Since $W'\neq 0$ in $\Pi_{-C,C}$ and is
real-valued on the real axis, one can write \be\log W'=\log
|W'|+i\Theta ,\label{lola}\ee where $\Theta$ is a harmonic function
in $\Pi_{-C, C}$, with $\Theta(x,0)=0$ for all $x\in\bdr$. Then \be
\frac{\partial U}{\partial x}=|W'|\cos\Theta,\quad \frac{\partial
V}{\partial x}=|W'|\sin\Theta.\label{fa}\ee Since
$W'\in\ha^1(\Pi_{-C,C})$ and $1/W'\in\ha^1(\Pi_{-C,C})$, it follows
that $\log |W'|\in\h^p(\Pi_{-C,C})$ for every $p\in (1,\infty)$ and
hence, by the M. Riesz Theorem, $\Theta\in\h^p(\Pi_{-C,C})$ for all
$p\in(1,\infty)$. Let $\te:\bdr\to\bdr$ be given, for almost every
$t\in\bdr$, by \be\te(t):=\lim_{(x,y)\toa
(t,C)}\Theta(x,y).\label{vte}\ee
 Let us also
denote, for all $t\in\bdr$, \be w(t)=u(t)+iv(t):=U(t,C)+i
V(t,C).\label{wewe} \ee Since $W'\in\ha^1(\Pi_{-C,C})$, it follows
from Lemma \ref{hac} that $w$ is locally absolutely continuous.
Moreover, taking into account (\ref{fa}), we deduce that \be
u'=|w'|\cos\te\quad\text{and}\quad v'=|w'|\sin\te\quad\text{almost
everywhere}.\label{ce1}\ee In particular, $\te$ gives the angle
between the tangent to the free boundary and the horizontal. The
condition (\ref{c6}) means that \be|w'|h(v)=1\quad\text{almost
everywhere}.\label{ce2}\ee It follows from (\ref{ce1}) and
(\ref{ce2}) that \be u'=\frac{\cos\te}{h(v)}\quad\text{and}\quad
v'=\frac{\sin\te}{h(v)}\quad\text{almost everywhere}.\label{ce3}\ee

Suppose now that $(\Om,\om)$ is a solution of problem {\bf (B)} as
in Theorem \ref{tlre}, and let $W$ be the corresponding solution of
problem {\bf (P)}. Let $t_1,t_2\in\bdr$ be such that
$\om(Z_j)=t_j+iC$ for $j=1,2$. It follows that $v(t)\in(c,d)$ for
all $t\in(t_1, t_2)$. As in \cite[Proof of Theorem 2.1 and Proof of
Theorem 2.3]{EV}, a simple bootstrap argument based on the second
part of Lemma \ref{lreg} and making use of (\ref{lola}), (\ref{fa}),
(\ref{vte}), (\ref{wewe}) and (\ref{ce3}) yields that $W\in
C^{n+1,\alpha}_{\textnormal{loc}}(\Pi_{0,C}\cup\{(t,C):t\in(t_1,t_2)\})$.
The required result is now immediate.
\end{proof}

\subsection{Proof of Theorem \ref{tgraph}}

\begin{proof}[Proof of Theorem \ref{tgraph}]

We make use of the notation and results in the proof of Theorem
\ref{tlre} concerning the corresponding solution $W$ of problem {\bf
(P)}. In particular, let $u,v$ be given by (\ref{wewe}). Then
clearly $\mcs=\{(u(t), v(t)):t\in\bdr\}$. Moreover, as we have seen,
$u$ and $v$ are locally absolutely continuous. It remains to prove
that (\ref{vepo}) holds.

Let $\mcn\subset\bdr$ be given by
\[\mcn:=\{t\in\bdr:h(v(t))=0\}=\{t\in\bdr:v(t)=d\}.\]  In the
terminology of \cite{ST,EV}, $\mcn$ is the set of stagnation points.
Obviously, $\mcn$ is a closed set and, in view of (\ref{ce2}), has
zero measure. By the proof of Theorem \ref{tlre}, $W\in
C^{2,\alpha}_{\textnormal{loc}}(\Pi_{0,C}
\cup\{(t,C):t\in\bdr\setminus\mcn\})$.

Note now that (\ref{ce2}) can be rewritten as \be -\log |w'|=\log
h(v)\quad\text{almost everywhere}.\label{becc}\ee Let
$Q:\Pi_{0,C}\to\bdr$ be given by \be Q(x,y)=  -\log |W'(x+iy)|-\log
h(V(x,y))\quad\text{for all }(x,y)\in \Pi_{0,C}.\label{iarb}\ee The
concavity of $\log h$ ensures that $Q$ is a subharmonic function in
$\Pi_{0,C}$. Also, it follows from (\ref{c4}) and (\ref{becc}) that,
in the notation of (\ref{Pois}),
\begin{align} Q&= \mathcal{P}_e[-\log |w'|]-\log
h(\mathcal{P}_o[v])\non\\&=\mathcal{P}_e[\log h(v)]-\log
h(\mathcal{P}_o[v]).\label{dern}\end{align} But since $v>0$, it
follows that \[ \mathcal{P}_e[v]>\mathcal{P}_o[v]>0.\] We deduce
from this and (\ref{nash}) by an obvious application of
 Jensen's Inequality \cite[Theorem 3.3, p.\ 62]{Ru} that
 \[ -\log h(\mathcal{P}_o[v])<-\log
h(\mathcal{P}_e[v])\leq \mathcal{P}_e[-\log h(v)].\] Now
(\ref{dern}) shows that $Q<0$ in $\Pi_{0,C}$. Since $Q$ is
subharmonic in $\Pi_{0,C}$, $Q\in C^1(\Pi_{0,C}\cup\{(x,C):
x\in\bdr\setminus\mcn\})$ and $Q=0$ in $\{(x,C):
x\in\bdr\setminus\mcn\}$,
 it follows from Hopf Boundary-Point Lemma that
\[\frac{\partial Q}{\partial y}(x,C)>0\quad\text{for all }x\in\bdr\setminus\mcn.\]
This means, upon using (\ref{iarb}) and the Cauchy-Riemann
equations, that \be
\te'(t)-\frac{h'(v(t))}{h^2(v(t))}\cos\te(t)>0\quad\text{for all
}t\in\bdr\setminus\mcn.\label{czam}\ee We aim to prove that \be
\cos\te(t)>0 \quad\text{for all
}t\in\bdr\setminus\mcn,\label{costa}\ee as this would yield
(\ref{vepo}).

 We argue by contradiction and assume that there exists
$f\in\bdr\setminus\mcn$ such that $\cos\te(f)\leq 0$. We distinguish
the following four cases:
\begin{itemize}
\item[(i)] $\mcn=\emptyset$,
\item[(ii)] $\mcn\neq\emptyset$ and $f\in(a,b)\subset\bdr\setminus\mcn$, where
$a,b\in\mcn$,
\item[(iii)]$\mcn\neq\emptyset$ and $f\in(a,\infty)\subset\bdr\setminus\mcn$, where
$a\in\mcn$,
\item[(iv)]$\mcn\neq\emptyset$ and $f\in(-\infty,b)\subset\bdr\setminus\mcn$, where
$b\in\mcn$.
\end{itemize}

We want to prove that a contradiction is reached in each of these
cases. We only give a full treatment of the case (i), and point out
the necessary modifications of the argument to deal with the
remaining cases.

As in \cite{EV2}, the proof in all the four cases here ultimately
rests on an application of the following classical theorem in the
global differential geometry of plane curves, see e.g.\ Amann
\cite[Theorem 24.15, p.\ 340]{Am} and the references therein, to a
suitably devised Jordan curve.

\begin{theorem}\label{TKRAS} Let $\sigma:[a,b]\to \bdc$ be a parametrization of a Jordan
curve, where $\sigma$ is a function of class $C^1$, with
$\sigma'(a)=\sigma'(b)$ and $|\sigma'|>0$ on $[a,b]$. Let
$\phi:[a,b]\to\bdr$ be a continuous function such that
\[\sigma'(t)=|\sigma'(t)|\exp\{i\phi(t)\}\quad\text{ for all }t\in[a,b].\]
Then $\phi(b)-\phi(a)$ equals either $2\pi$ or $-2\pi$.
\end{theorem}

 We also use the following lemmas \cite[Proof of Lemma 4.16 and Proof of Lemma 4.17]{ST}.

\begin{lemma}\label{LISO} Let $(a, b) \subset \bdr \setminus \mcn$
and suppose that  {\rm (\ref{czam})} holds on $(a, b)$. Let $e\in
(a,b)$. If $\ell_1, \ell_2 \in\ZZ$ are such that
$\ell_1\pi+\pi/2\leq\te(e)\leq\ell_2\pi+\pi/2$, then
\[ \te(t)
> \ell_1\pi +\pi/2\text{ for } t\in (e,b)\quad\text{ and }\quad\te(t)
<\ell_2\pi +\pi/2\text{ for } t\in (a,e).\]
\end{lemma}

\begin{lemma}\label{LBOUND} Let $(a, b) \subset \bdr \setminus
\mcn$, where $a\in\mcn\cup\{-\infty\}$ and
$b\in\mcn\cup\{+\infty\}$, and suppose that {\rm (\ref{czam})} holds
on $(a, b)$. Let $f\in (a,b)$ be such that $\cos \te(f)\leq 0$. If
$m\in\ZZ$ is such that
\[2m\pi+\pi/2\leq \te(f)\leq 2m\pi+3\pi/2,\] then there exists
$g_1,g_2\in (a,b)$ with $g_1\leq f\leq g_2$ such that $\te$ is
strictly increasing on $[g_1,g_2]$ and
\[ \te(g_1)=2m\pi+\pi/2\quad\text{ and }\quad \te(g_2)=2m\pi+3\pi/2.\]
\end{lemma}

\begin{proof}[Proof of Lemma \ref{LBOUND}]
We prove only the existence of $g_2$ with the required properties,
since the proof for $g_1$ is entirely similar.

 When $b\in\mcn$, the existence of $g_2$ follows by the argument
 in \cite[Proof of Lemma 4.17]{ST}.

Suppose now that $b=+\infty$. Note from (\ref{czam}) that $\te$ is
strictly increasing on each interval in which $\cos\te\leq 0$. If
such $g_2$ does not exist, then obviously $\cos\te <0$ on $(f,
\infty)$, and this implies that $u'<0$ on $(f,\infty)$, which
contradicts the fact that $\lim_{t\to\infty}u(t)=+\infty$. This
proves the existence of $g_2$.
\end{proof}

In what follows we deal with the case (i), thus assuming that
$\mcn=\emptyset$ and there exists $f\in\bdr$ such that
$\cos\te(f)\leq 0$. It follows that \be \te(f)\in[2m\pi+\pi/2,
2m\pi+3\pi/2]\quad\text{ for some }m\in\mathbb{Z}.\label{ttz}\ee
 Let
$g_1$, $g_2$ be given by Lemma \ref{LBOUND}, so that \be
\te(g_1)=2m\pi+\pi/2,\quad\te(g_2)=2m\pi+3\pi/2.\label{gud}\ee Let
$M_1,M_2\in\bdr$ be such that
\[M_1<\min\{u(s):s\in[g_1,g_2]\},\quad M_2>\max\{u(s):s\in[g_1, g_2]\}.\]
Let $p_1, p_2$ be such that \[ u(p_1)=M_1,\quad u(p_2)=M_2,\quad
M_1<u(s)<M_2\text{ for all }s\in(p_1, p_2).\] Then
\[ u'(p_1)\geq 0,\quad u'(p_2)\geq 0,\]
so that there exists $m_1, m_2\in\ZZ$ such that
\be\te(p_1)\in[2m_1\pi-\pi/2,2m_1\pi+\pi/2],\quad
\te(p_2)\in[2m_2\pi-\pi/2,2m_2\pi+\pi/2].\label{tepe}\ee It follows
from (\ref{gud}), (\ref{tepe}) and Lemma \ref{LISO} that \be m_1\leq
m,\quad m_2\geq m+1.\label{mud} \ee

 Let $\tilde w:[p_1,p_2]\to\bdc$
be the restriction of $w$ to $[p_1,p_2]$. It is obvious that, for
some
  $q_1,\,q_2\in\bdr$ with $q_1<p_1$, $p_2<q_2$, one can construct a function
  $\hat w:[q_1, q_2]\to\bdc$, where
\[\hat w(q):=\hat u(q)+i\hat v(q)\quad\text{for all }q\in[q_1,q_2],\]
 such that $\hat w$ is an extension of $\tilde w$, and
it has the following additional properties: \bese\label{multe}
\begin{align}\hat u, \hat v:[q_1,q_2]&\to\bdr\quad \text{ are of class }\,C^1,\\\hat u'(q)^2+\hat v'(q)^2
  &>0\quad\text{ for all }q\in[q_1,q_2],\\
   \hat u'&\geq 0\quad\text{ on }[q_1,p_1]\cup [p_2,q_2],\\\hat u'(q_1)= \hat u'(q_2)&=1,\\
  \hat v'(q_1)= \hat v'(q_2)&=0.
\end{align}\ese

Let $E$ be such that
\[E>\max\{\hat v(q):q\in[q_1, q_2]\}.\]
Let $\mathcal{A}_1$ be the semicircle having as diameter the segment
joining the points $(\hat u(q_1), E)$ and $(\hat u(q_1), \hat
v(q_1))$, and situated to the left of this segment. Let
$\mathcal{A}_2$ be the semicircle having as diameter the segment
joining the points $(\hat u(q_2),\hat v(q_2))$ and $(\hat u(q_2),
E)$, and situated to the right of this segment. Let $r_1,\,r_2$ with
$r_1<q_1$, $r_2>q_2$ and
 consider a $C^1$ function $w_*:[r_1, r_2]\to\bdc$ which is an extension
  of $\hat w$, such that\bese\label{mult1}
 \begin{align}w_*(r_1)=\hat u(q_1)+iE&,\quad w_*(r_2)=\hat u(q_2)+iE,
 \\w_*'(r_1)&=w_*'(r_2)=-1+i0,\\
 |w_*'(r)|&>0\quad\text{ for all }
r\in[r_1,r_2],\end{align}\ese
 and \[\begin{aligned} w_*|_{[r_1,q_1]}\text{ is an injective
parametrization of  $\mathcal{A}_1$},\\
  w_*|_{[q_2,r_2]}\text{ is an injective parametrization of $\mathcal{A}_2$}.\end{aligned}\]
Let $\tilde r_1<r_1$ be such that $w_*$ extends to $[\tilde r_1,
r_1]$ as a $C^1$ function such that \[w_*(\tilde r_1)=\hat
u(q_2)+iE,\]
\[w_*'(r)=-1+i0\quad\text{ for all } r\in[\tilde r_1, r_1].\]
and let also $\tilde r_2:=r_2$.

It is very easy to prove that $w_*:[\tilde r_1, \tilde r_2]\to\bdc$
constructed above provides a parametrization of a Jordan curve with
a continuously varying tangent. Let us write
\[ w_*'(r)=|w_*'(r)|\exp\{i\te_*(r)\}\quad\text{ for all } r\in[\tilde r_1, \tilde r_2],\]
where $\te_*:[\tilde r_1, \tilde r_2]\to\bdr$ is a continuous
function which extends $\te:[p_1,p_2]\to\bdr$. It follows from
(\ref{tepe}) and (\ref{multe}) that
\[\te_*(q_1)=2m_1\pi,\qquad\te_*(q_2)=2m_2\pi.\]
Using (\ref{mult1}) we deduce that
\[\te_*(\tilde r_1)=2m_1\pi-\pi,\qquad\te_*(\tilde r_2)=2m_2\pi+\pi.\]
Therefore \be\te_*( \tilde r_2)-\te_*(\tilde
r_1)=2(m_2-m_1)\pi+2\pi,\label{difg}\ee where, by (\ref{mud}),
\be\label{difg1}m_2-m_1\geq 1.\ee But the validity of (\ref{difg})
with (\ref{difg1}) is in contradiction to Theorem \ref{TKRAS}. This
proves the required result (\ref{costa}) in the case (i).

The case (ii) can be dealt with by an argument which is entirely
similar to that in \cite[Proof of Theorem 3.3]{EV2}. The cases (iii)
and (iv) can be treated using a geometric construction which
combines elements of those used for the cases (i) and (ii). This
completes the proof of Theorem \ref{tgraph}.
\end{proof}

{\bf Acknowledgement.} I am very grateful to Professor J.F. Toland
for many valuable discussions on the theory of water waves. This
work was supported by the EPSRC grant no.\ EP/D505402/1.

\end{document}